\documentclass[12pt,twoside]{article}

\date{}

\begin{document}

\centerline{}

\centerline {\Large{\bf A Note on Hypervector Spaces}}
\newcommand{\mvec}[1]{\mbox{\bfseries\itshape #1}}
\centerline{}
\centerline{\bf {Sanjay Roy, T.K. Samanta }}

\centerline{Department of Mathematics}
 \centerline{South Bantra Ramkrishna Institution , West Bengal, India. }
\centerline{e-mail: sanjaypuremath@gmail.com}

\centerline{Department of Mathematics, Uluberia College, West Bengal, India.}
\centerline{e-mail: mumpu$_{-}$tapas5@yahoo.co.in}

\newtheorem{Theorem}{\quad Theorem}[section]

\newtheorem{definition}[Theorem]{\quad Definition}

\newtheorem{theorem}[Theorem]{\quad Theorem}

\newtheorem{remark}[Theorem]{\quad Remark}

\newtheorem{corollary}[Theorem]{\quad Corollary}

\newtheorem{note}[Theorem]{\quad Note}

\newtheorem{lemma}[Theorem]{\quad Lemma}

\newtheorem{example}[Theorem]{\quad Example}

\newtheorem{result}[Theorem]{\quad Result}

\newtheorem{proposition}[Theorem]{\quad Proposition}

\begin{abstract}
 \textbf{\emph{The main aim of this paper is to generalize the concept of vector space by the hyperstructure. We generalize some definitions such as hypersubspaces, linear combination, Hamel basis, linearly dependence and linearly independence. A few important results like deletion theorem, extension theorem, dimension theorem have been established in this hypervector space.}}
\end{abstract}
{\bf Keywords:}  \emph{Hyperoperation, Hyperfield, Hypervector spaces, linear dependent, linear independent.}\\
\textbf{2010 mathematics subject classification:} 14L17, 20N20.

\section{\textbf{Introduction}}
 The concept of hyperstructure was first introduced by Marty \cite{Marty1} in 1934 at the 8th congress of scandinavian Mathematicians and then he established the definition of hypergroup \cite{Marty2} in 1935 to analyse its properties and applied them to groups of rational algebraic functions. Also he was motivated to introduce this structure to study several problems of the non-commutative algebra. Then several researchers have been working on this new field of modern algebra and developed it. M. Krasner \cite{Krasner}, a great researcher in this area, introduced the notions of hyperring and hyperfield to use it as a technical tool in the study of the approximation of valued fields. Later on it has been developed and generalized by other researchers. Then the notion of the hypervector spaces was introduced by M. Scafati Tallini \cite{Rivista} in 1988.
\\\\In the definition\cite{M.s.Tallini} of hypervector spaces, M.Scafati Tallini has considered the field as a usual field. In this paper, we have generalized the definition of hypervector space by considering the field as a hyperfield and considering the  multiplication structure of a vector by a scalar as hyperstructure like M. Scafati Tallini. We again call it a hypervector space. Then we have established a few basic properties in this hypervector space and thereafter the notions of linear combinations, linearly dependence, linearly independence, Hamel basis, etc. are introduced and several important properties like deletion theorem, extension theorem etc. are developed.

\section{\textbf{Preliminaries}}
 We quote some definitions and proofs of a few results which will be needed in the sequel.

\begin{definition}\cite{Nakassis} A \textbf{hyperoperation} over a non-empty set X is a mapping of X$\times $X into the set of all non-empty subsets of X.
\end{definition}

\begin{definition}\cite{Nakassis}
A non-empty set X with exactly one hyperoperation `$\#$' is a \textbf{hypergroupoid}. \\
Let $(X\, ,\, \#)$ be a hypergroupoid. For every point x $\in$ X and every non-empty subset A of X, we define x $\#$ $A=\bigcup_{a\in A}\{ x \;\#\; a\}$.                                                                                                                             \end{definition}

\begin{definition}\cite{Nakassis}
A hypergroupoid $(X \,,\,\#)$ is called a \textbf{semihypergroup} if $x\#(y\# z)$ = $(x\#y)\# z$  for all x, y, z $\in X$.
\end{definition}

\begin{definition}\cite{Nakassis} A hypergroupoid $(X \,,\,\#)$ is called a \textbf{hypergroup} if\\
( $i$ )\hspace{0.5cm}   $x\#(y\#z)=(x\#y)\#z$\\
( $ii$ )\hspace{0.4cm}  $\exists$ $0$ $\in$ X such that for every a $\in$ X, there is unique element b $\in$ X for which $0\in a\#b$ and $0\in b\# a$ . Here $b$ is denoted by $- a$.\\
( $iii$ ) \hspace{0.3cm}For all a, b, c $\in X$ if $a\in b \# c$, then $b \in a\#(- c)$.
\end{definition}

\begin{result}\cite{Nakassis}
$:$ In a hypergroup $(X \,,\, \#)$, $-(- a) = a$, $\forall  a \in X$.
\end{result}
\textbf{Proof} Since $0 \in a\#(-a)$ and $0\in (-a)\#a,\;\forall a\in X$\\
i.e $0 \in(-a)\# a$ and $0 \in a\#(-a),\; \forall a \in X$.\\
Hence $-(-a)=a$,\,\,$\forall a\in$ X.

\begin{result}\cite{Nakassis} $0\# a=\{a\}$, $\forall a \in X$, if $(X\,,\,\#)$ is a commutative hypergroup.
\end{result}
\textbf{Proof} Let $a \in X$, then $0 \in a \#(-a)$\\
$\Rightarrow\; a \in 0 \#(-(-a))$    $[$by definition 2.4 $]$\\
$\Rightarrow\; a \in 0\# a $ \\
we now show that $0\#a=\{a\}$\\
Let b $\in 0 \# a$\\
$\Rightarrow\, 0 \in b \#(-a)=(-a)\# b$  \\
$\Rightarrow\, b=-(-a)$\\
$\Rightarrow\; b=a.$\\
This completes the proof.

\begin{result}\cite{Nakassis} In a commutative hypergroup $($X , $\#)$, $0$ is unique.
\end{result}
\textbf{Proof} Let there be another element e $\in$ X such that\\
 e $\in a \;\#\;(-a)$ , e $\in (-a)\;\#\;a,\;\forall\, a\in X$.\\
$\Rightarrow\,a \#\, e=\{a\},$ $\forall\, a\in$X.\\
$\Rightarrow$  $0 \,\#\, e=\{0\}$. Again $0 \,\# \,e=\{e\}$. \\
Hence e=0.\\
This completes the proof.

\begin{note} In a hypergroup, if the element $0$ is unique , then $0$ is called the zero element of the hypergroup and b is called the additive inverse of a if $0 \in a\#b$ and $0\in b\# a$.
\end{note}

\section{\textbf{Hypervector Space}}
\begin{definition}A \textbf{hyperring} is a non-empty set endowed with a hyperaddition `\,$\oplus$' and a multiplication `.' such that $(X\,,\,\oplus)$ is a commutative hypergroup and $(X\,,.)$ is a semigroup and the multiplication is distributive with respect of the hyperaddition, both from the left and from the right side and $a.0\,=\,0.a\,=\,0$  $\forall\, a\in X$ , where $0$ is the zero element of the hyperring.
\end{definition}
\begin{definition}
A \textbf{hyperfield} is a non-empty set X endowed with a hyperaddition `\,$\oplus$' and a multiplication `.' such that\\
$(\,i\,)$\hspace{0.5cm}$(X\,,\;\oplus\;,\,.)$ is a hyperring.\\
$(\,ii\,)$\hspace{0.4cm}$\exists$ an element $1 \in X$ , called the identity element such that $a.1=a$ , $\forall\, a\in X.$ \\
$(\,iii\,)$\hspace{0.3cm}For each non zero element $a \in X$ , $\exists$ an element $a^{-1}\in X$ such that $a.a^{-1}=1$.\\
$(\,iv\,)$\hspace{0.4cm}$a.b=b.a$ , $\forall\, a , b\in X$.
\end{definition}
\begin{definition}Let $(F\,,\oplus\,,.)$ be a hyperfield and $(V\,,\,\#)$ be an additive commutative hypergroup . Then V is said to be a \textbf{hypervector space} over the hyperfield F if there exist a hyperoperation $\ast:\, F \,\times\, V \,\rightarrow\;P^{\ast}\,(\,V\,)$ such that\\
$(\,i\,)$\hspace{0.5cm}a$\ast(\alpha\;\#\;\beta)\subseteq a\ast\alpha\;\,\#\;\,a\ast\beta$ , \hspace{1cm}$\forall\, a\in F\;and\;\forall \alpha\, , \beta \in$ V.\\
$(\,ii\,)$\hspace{0.4cm}$(a\,\oplus\,b)\ast\alpha \subseteq a\ast\alpha\;\#\;b\ast\alpha$  , \hspace{1cm} $\forall\,a , b\in F\;and\;\forall\,\alpha\in
$ V.\\
$(\,iii\,)$\hspace{0.3cm}$(a\,.\,b)\ast\alpha=a\ast(b\ast\alpha)$  , \hspace{2cm} $\forall\, a , b \in F$ and $\forall \alpha\in$ V.\\
$(\,iv\,)$\hspace{0.4cm}$1_{F}\ast\alpha=\{\alpha\}\; and\;0\ast\alpha=\{\theta\}$  ,\hspace{1cm}  $\forall \alpha\in$ V where $1_{F}$ is the identity element of F , 0 is the zero element of F and $\theta$ is the zero vector of V and $P^{\ast}(\,V\,)$ is the set of all non-empty subsets of V.\\\\
A hypervector space is called \textbf{strongly right distributive hypervector space} $($ respectively , \textbf{strongly left distributive hypervector space} $)$ , if eqality holds in $(\,i\,)\;( $ respectively , $in\;(\,ii\,)\;)$.\\
A hypervector space is called a \textbf{good hypervector space} if equality holds in both $(\,i\,)\;and\;(\,ii\,)$.
\end{definition}
\begin{remark}By a hypervector space V, we mean a hypervector space $(V\,,\#\,,\ast)$ and by a hyperfield F, we mean a hyperfield $(F\,,\oplus\,,.)$.
\end{remark}
\begin{remark}Let V be a hypervector space over a hyperfield F. Let $a\;,\;b \in$ F and $\alpha\;,\beta\in$ V. then by $ a\ast\alpha\;\#\;b\ast\beta$, we mean $(a\ast\alpha)\;\#\;(b\ast\beta)$.
\end{remark}
\begin{example}Let $(F\,,\,\oplus\,,.)$ be an hyperfield and $V = F \times$ F. Let us define a hyperoperation `$\#$' on V as follows\\
$(a_{1}\,,\,a_{2})\;\#\;(b_{1}\,,\,b_{2})=(a_{1}\oplus b_{1}\;,\;a_{2}\oplus b_{2})$=$\{ (x ,y) : x\in a_{1}\oplus b_{1}$ and y $\in a_{2}\oplus b_{2}\}.$\\
Then we prove that $(V\;,\;\#)$ is an additive commutative hypergroup. Now we define a scalar multiplication $\ast:\;F\;\times\;V\;\rightarrow\;P^{\ast}(\,V\,)$ by \\
a$\ast((a_{1}\,,\,a_{2}))=\{(a.a_{1}\,,\,a.a_{2})\}$ , where a $\in$ F and $(a_{1}\,,\,a_{2})\in$ V. Then we easily verify that \\
$(\,i\,)\;a\ast((a_{1}\,,\,a_{2})\;\#\;(b_{1}\,,\,b_{2}))=(a\ast(a_{1}\,,\,a_{2}))\;\#\;(a\ast(b_{1}\,,\,b_{2})),$\\
$(\,ii\,)\;(a\oplus b)\ast(a_{1}\,,\,a_{2})=(a\ast(a_{1}\,,\,a_{2}))\;\#\;(b\ast(a_{1}\,,\,a_{2})),$\\
$(\,iii\,)\;(a.b)\ast(a_{1}\,,\,a_{2})=a\ast(b\ast(a_{1}\,,\,a_{2})),$\\
$(\,iv\,)\;1_{F}\ast(a_{1}\,,\,a_{2})=(a_{1}\,,\,a_{2})$ and $0\ast(a_{1}\,,\,a_{2})=(0\,,\,0)=\theta,$
for all $a , b \in\;F$ and for all $(a_{1}\,,\,a_{2})\;,\;(b_{1}\,,\,b_{2})\in$ V.
\end{example}
\begin{result}Let $(V\,,\;\#\;,\,\ast)$ be a hypervector space over a hyperfield $(F\,,\,\oplus\,,\,.)$. Then\\
$(\,i\,)$\hspace{0.5cm}k $\ast\theta=\{\theta\}$ , $\forall\;k\in$ F , $\theta$ being the zero vector of V.\\
$(\,ii\,)$\hspace{0.4cm}Let k $\in$ F and $\alpha\in$ V be such that k $\ast\alpha=\{\theta\}$ , then either k=0 ,or $\alpha=\theta.$\\
$(\,iii\,)$\hspace{0.4cm}$-\alpha\,\in\,(-1_{F})\ast\alpha\;,\;\forall\;\alpha\in$ V , $1_{F}$ being the identity element of F.
\end{result}
\textbf{Proof}:$(i)\; k\ast\theta=k\ast(0\ast\theta)$ , $[$ by axiom $(\,iv\,)$ , we have 0$\ast\alpha=\theta\;,\;\forall\;\alpha \in$ V $]$\\
\smallskip\hspace{2.7cm}=$(k\,.\,0)\ast\theta\,=\,0\ast\theta\,=\,\theta.$\\
\textbf{Proof}:$(ii)\;$ Let $k \in$ F and $\alpha\;\in$ V be such that k $\ast
\;\alpha\;=\;\{\theta\}.$\\
If $k=0$ , then $0\ast\alpha=\theta.$\\
If $k\neq 0$ , then $k^{-1}\in$ F.\\
Therefore $k\ast\alpha=\theta\;\Rightarrow\;k^{-1}\ast(k\ast\alpha)=k^{-1}\ast\theta$\\
\smallskip\hspace{3.6cm}$\Rightarrow \;(k^{-1}\,.\,k)\ast\alpha=\theta\,\Rightarrow\;1_{F}\ast\alpha=\theta\,\Rightarrow\;\alpha=\theta$. \\
This completes the proof.\\
\textbf{Proof}:$(iii)\;$  Let $\alpha\in V$ , then $(1_{F}\oplus(-1_{F}))\ast\alpha\,\subseteq\,1_{F}\ast\alpha\;\#\;(-1_{F})\ast\alpha$\\
 \smallskip\hspace{8.6cm}= $\alpha\;\#\;(-1_{F})\ast\alpha.$\\
 Since 0 $\in\;1_{F}\oplus(-1_{F})\;\Rightarrow\;\theta=0\ast\alpha\in(1_{F}\oplus(-1_{F}))\ast\alpha$\\
 \smallskip\hspace{4.2cm}$\Rightarrow\;\theta\in\alpha\;\#\;(-1_{F})\ast\alpha$.\\
 Therefore  $-\alpha\,\in\,(-1_{F})\ast\alpha$  ,  $\forall\;\alpha\;\in$ V.

\section{\textbf{HypersubSpaces}}
\begin{definition} A subset W of a hypervector space V over a hyperfield F is called a \textbf{hypersubspace} of V if W is a hypervector space over F with the hyperoperations of addition and the scalar multiplication defined on V.\\
  Therefore a subset W of a hypervector space V is a hypersubspace of V if and only if the following four properties hold.\\
  $(\,i\,)$\hspace{0.5cm}   $\alpha\;\#\;\beta\subseteq W,\;\; \forall\,\alpha ,\beta\,\in\, W$.\\
  $(\,ii\,)$\hspace{0.4cm}  $a\ast\alpha\,\subseteq\,W,\;\;\forall\,\alpha\in W\;and\;\forall\;a\;\in$ F.\\
  $(\,iii\,)$ \hspace{0.3cm}W has a zero vector.\\
  $(\,iv\,)$\hspace{0.4cm} each vector of W has an additive inverse.
  \end{definition}
  \begin{theorem} Let V be a hypervector space and W is a subset of V . Then W is a hypersubspace of V if and only if the following three conditions hold$:$\\
  $(\,i\,)$ W is non-empty.\\
  $(\,ii\,)\;\alpha\;\#\;\beta\,\subseteq$ W  , $\forall\,\alpha\ , \beta\in$ W.\\
  $(\,iii\,)\;a\ast\alpha\,\subseteq$ W , $\forall\,a\in$ F and $\forall\,\alpha\in$ W.
  \end{theorem}
  \textbf{Proof}: If W is a hypersubspace of V , then obviously the conditions $(i)\;,(ii)\;and\;(iii)$ hold.\\
  Conversely , let W be a subset of V such that W satisfies the three conditions $(i)\;,(ii)\;and\;(iii)$.\\
  To proof that W is a hypersubspace of V. it is enough to prove that \\
   $(1)$ W has a zero vector. $(2)$Each vector in W has an additive inverse.\\
  Since W is non-empty , let $\alpha\in$ W.\\
  Now $0\in$ F , therefore by the condition $(iii)$ we get $0\ast\alpha\;\subseteq$ W $\Rightarrow\;\theta\in$ W.\\
  Therefore W has a zero vector.\\
  Again, since $-1_{F}\in$ F,\\
 therefore $(-1_{F})\ast\alpha\,\subseteq$ W $\Rightarrow\;-\alpha\in$ W.\\
 Hence each vector in W has an additive inverse.

 \begin{theorem}W be a hypersubspace of a hypervector space V if and only if $(\,i\,)$ W is non-empty. $(\,ii\,)\;a\ast\alpha\;\#\;b\ast\beta\,\subseteq$ W , $\forall$ a, b $\in$ F and $\forall\;\alpha, \beta\in$ W.
 \end{theorem}
 \textbf{Proof}:  If W is a hypersubspace of V, then obviously W satisfies the conditions $(i)\;and\;(ii)$.\\
 Conversely , let W satisfy the conditions $(i)$ and $(ii).$\\
 Since $1_{F}\in$ F , let $\alpha ,\beta\in$ W , then by $(ii)$\\
 $1_{F}\ast\alpha\;\#\;1_{F}\ast\beta\,\subseteq$ W $\Rightarrow\;\alpha\;\#\;\beta\,\subseteq$ W.
 Let $a$ $\in$ F and $\alpha\,,\,\beta\in$ V , since $0\in$ F , therefore by $(ii)$\\
 $a \ast\alpha\;\#\;0\ast\beta\,\subseteq W\;\Rightarrow\;a\ast\alpha\;\#\;\theta\,\subseteq$ W $\Rightarrow\;a\ast\alpha\,\subseteq$ W  $[\;Because \;(W\,,\;\#)$ is a commutative hypergroup $]$\\
  Hence W is a hypersubspace of V.

\begin{example}$:$ Let $(F\,,\,\oplus\,,\,.)$ be a hyperfield and V = F $\times$ F.\\
Then $(V\,,\;\#\;,\,\ast)$ is a hypervector space, where the hyperoperations $'\#'\;and\;'\ast'$ are defined by\\
$(a_{1},a_{2})\;\#\;(b_{1},b_{2})=(a_{1}\oplus b_{1}\;,\;a_{2}\oplus b_{2})$=$\{ (x ,y) : x\in a_{1}\oplus b_{1}$ and y $\in a_{2}\oplus b_{2}\}$ \\
and $ a\ast(a_{1} , a_{2})=\{(a.a_{1} , a.a_{2})\}$ ,  $\forall\;(a_{1} , a_{2})\;,\;(b_{1} , b_{2})\in V$ and $\forall\;a\in$ F.\\
Let W = F $\times\;\{0\}\,\subseteq$ V\\
We now show that W is a hypersubspace of V,\\
since $\theta=(0 , 0)\in$ W.\\
Now let $\alpha=(a_{1} , 0)\; ,\;\beta=(b_{1} , 0)\in W$ and $a , b\,\in\, F$.\\
Then a$\ast\alpha\;\#\;b\ast\beta$\\
 = a$\ast(a_{1} , 0)\;\#\;b\ast(b_{1} , 0)$\\
 = $\{(a.a_{1} , a.0)\}\;\#\;\{(b.b_{1} , b.0)\}$\\
 = $\{(a.a_{1} , 0)\}\;\#\;\{(b.b_{1} , 0)\}$\\
 = $(a.a_{1}\oplus b.b_{1}\;,\;0\oplus 0)$\\
 = $(a.a_{1} \oplus b.b_{1}\;,\;0)\subseteq$ W  $[Because\;a.a_{1}\;\oplus\;b.b_{1}\subseteq F].$
 \end{example}
\begin{theorem} The intersection of two hypersubspaces of a hypervector space V over a hyperfield F is again a hypersubspace of V.
\end{theorem}
\textbf{Proof}: Obvious.

\begin{theorem}The intersection of any family of hypersubspaces of a hypervector space V over a hyperfield F is again a hypersubspace of V.
\end{theorem}
\textbf{Proof}: Obvious.

  In the above example, we take $W_{1}$ = F $\times \{0\}\,\subseteq V$ and $W_{2}=\{0\}\times F\subseteq V.$\\
Then by the same procedure of the above example , we can prove that $W_{1}$ and $W_{2}$ are hypersubspaces of V.\\
Let $(a , 0)\in W_{1}$ and $(0 , b)\in W_{2}.$\\
Then $(a , 0)\;\#\;(0 , b)=(a\oplus0\;,\;0\oplus b)=(a\;,\;b).$\\
But $\{(a , b)\}$ is not a subset of $W_{1}\cup W_{2}$\\
Therefore $W_{1}\cup W_{2}$ is not a hypersubspace of V.

\begin{note}$:$ The union of two hypersubspaces of a hypervector space V is not necessarily a hypersubspace of V.
\end{note}
\begin{theorem} Let $W_{1}$ and $W_{2}$ be two hyperspaces of a hypervector space V. Then we prove that $W_{1}\;\#\;W_{2}=\cup\{\alpha\;\#\;\beta\;,\;\alpha\in W_{1}\;,\;\beta\;\in\;W_{2}\}$ is a hypersubspace of V.
\end{theorem}
\textbf{Proof}: Since $\theta\in W_{1}$ and $\theta\in W_{2}.$\\
Then $\{\theta\;\#\;\theta\}\,\subseteq W_{1}\;\#\;W_{2}\;\Rightarrow\;\{\theta\}\,\subseteq\,W_{1}\;\#\;W_{2}\;\Rightarrow\;\theta\in W_{1}\;\#\;W_{2},$\\
therefore $W_{1}\;\#\;W_{2}$ is non-empty.\\
Let $\alpha , \beta\in W_{1}\;\#\;W_{2}$ , Then $\exists\, \alpha_{1} , \alpha_{2}\in W_{1}$ and $\beta_{1} , \beta_{2}\in W_{2}$ such that $\alpha\in \alpha_{1}\;\#\;\beta_{1}$ and $\beta\in \alpha_{2}\;\#\;\beta_{2}$. Let $a$ , $b$ $\in$ F.\\
Now $a\ast\alpha\;\# \;b\ast\beta\,\subseteq\,a\ast(\alpha_{1}\;\#\;\beta_{1})\;\#\;b\ast(\alpha_{2}\;\#\;\beta_{2})$\\
\smallskip\hspace{3cm}$\subseteq\;(a\ast \alpha_{1}\;\#\;a\ast \beta_{1})\;\#\;(b\ast \alpha_{2}\;\#\;b\ast \beta_{2})$\\
\smallskip\hspace{3cm}=$a\ast \alpha_{1}\;\#\;(a\ast \beta_{1}\;\#\;b\ast \alpha_{2})\;\#\;b\ast \beta_{2}$\\
\smallskip\hspace{3cm}=$a\ast \alpha_{1}\;\#\;(b\ast \alpha_{2}\;\#\;a\ast \beta_{1})\;\#\;b\ast \beta_{2}$\\
\smallskip\hspace{3cm}=$(a\ast \alpha_{1}\;\#\;b\ast \alpha_{2})\;\#\;(a\ast \beta_{1}\;\#\;b\ast \beta_{2})$\\
\smallskip\hspace{3cm}$\subseteq\,W_{1}\;\#\;W_{2}.$\\
Hence $W_{1}\;\#\;W_{2}$ is a hypersubspace of V.

\begin{definition}If $W_{1}$ and $W_{2}$ are two hypersubspaces of a hypervector space V, then the hypersubspace $W_{1}\;\#\;W_{2}$ is called the \textbf{hyperlinear sum} or \textbf{linear sum} of the hyperspaces $W_{1}$ and $W_{2}$.\\
     If $W_{1}\cap W_{2}=\{\theta\}$  then $W_{1}\;\#\;W_{2}$ is called the \textbf{direct sum} of the hypersubspaces $W_{1}$ and $W_{2}$.
\end{definition}
\begin{theorem}
The hypersubspace $W_{1}\;\#\;W_{2}$ is the smallest hypersubspace of V containing the hypersubspaces $W_{1}$ and $W_{2}$.
\end{theorem}
\textbf{Proof}:  Let W be a hypersubspace of V such that $W_{1}\;\subseteq$ W and $W_{2}\;\subseteq$ W. Let $\gamma\in W_{1}\;\#\;W_{2}$ , then $\exists$ $\alpha\in\;W_{1}$ and $\beta\in\;W_{2}$ such that $\gamma\in\;\alpha\;\#\;\beta$\\
Since $W_{1}\,\subseteq W$ and $W_{2}\,\subseteq$ W. Therefore $\alpha\; , \beta \in$ W. Again since W is a hypersubspace of V.\\
Therefore $\alpha\;\#\;\beta\;\subseteq\;W\;\Rightarrow\;\gamma\in$ W.\\
Hence $W_{1}\;\#\;W_{2}\,\subseteq$ W\\
This completes the proof.\\

If V is a strongly left distributive hypervector space over a hyperfield F, then it can be easily prove that W = $\cup\{a\ast\alpha,\; a\in F\}$ forms a hypersubspace of V, where $\alpha
\in$ V. This hypersubspace is said to be generated by the vector $\alpha$ and $\alpha$ is said to be a \textbf{generator} of the hypersubspace W. This hypersubspace is usually denoted by HL$(\alpha)$.\\
\smallskip\hspace{0.2cm}Then from the previous theorem it can be easily proved that if $\alpha , \beta\in$ V , then the set W =$ \cup\{a\ast\alpha\;\#\;b\ast\beta,\; a,\,b \in$ F$\}$ is a hypersubspace of V. This hypersubspace is called the \textbf{hyperlinear span} of the vectors $\alpha$ and $\beta$ , it is usually denoted by HL$(\alpha , \beta)$.
\begin{theorem} Let V be a strongly left distributive hypervector space over the hyperfield F and $\alpha_{1} , \alpha_{2} ,.\;.\;.\;,\;\alpha_{n}\in$ V. Then\\
W = $\cup\{a_{1}\ast\alpha_{1}\;\#\;a_{2}\ast\alpha_{2}\;\#\;.\;.\;.\;\#\;a_{n}\ast\alpha_{n}\;:\;a_{1}\;,a_{2}\;,\;.\;.\;.\;,\;a_{n}\;\in\;F\}$ is a hypersubspace of V. In fact W is the smallest hypersubspace of V containing $\alpha_{1} , \alpha_{2} ,\;.\;.\;.\;, \alpha_{n}$.
\end{theorem}
\textbf{Proof}: Since F is non-empty, therefore W is non-empty. Let $w_{1} , w_{2}\;\in$ W. Then $\exists\;a_{1} , a_{2} , .\;.\;.\;, a_{n}\;,\;b_{1} , b_{2} ,\;.\;.\;.\;, b_{n}\in F$ such that \\
$w_{1}\in a_{1}\ast\alpha_{1}\;\#\;a_{2}\ast\alpha_{2}\;\#\;.\;.\;.\;\#\;a_{n}\ast\alpha_{n}$ and\\
$w_{2}\in b_{1}\ast\alpha_{1}\;\#\;b_{2}\ast\alpha_{2}\;\#\;.\;.\;.\;\#\;b_{n}\ast\alpha_{n},$\\
therefore $w_{1}\;\#\;w_{2}\,\subseteq\,(a_{1}\ast\alpha_{1}\;\#\;a_{2}\ast\alpha_{2}\;\#\;.\;.\;.\;\#\;a_{n}\ast\alpha_{n})\;\#\;
(b_{1}\ast\alpha_{1}\;\#\;b_{2}\ast\alpha_{2}\;\#\;.\;.\;.\;\#\;b_{n}\ast\alpha_{n})$\\
or $w_{1}\;\# w_{2}\,\subseteq\,(a_{1}\ast\alpha_{1}\;\#\;b_{1}\ast\alpha_{1})\;\#(a_{2}\ast\alpha_{2}\;\#\;b_{2}\ast\alpha_{2})\;\#\;.\;.\;.\;\#
(a_{n}\ast\alpha_{n}\;\#\;b_{n}\ast\alpha_{n})$ \\
or $w_{1}\;\# w_{2}\;\subseteq\;(a_{1}\oplus b_{1})\ast\alpha_{1}\;\#\;(a_{2}\oplus b_{2})\ast\alpha_{2}\;\#\;.\;.\;.\;\#(a_{n}\oplus b_{n})\ast\alpha_{n}.$\\
Therefore $w_{1}\;\#\;w_{2}\;\subseteq$ W\\
Next let $w \in$ W and $a \in$ F. Then $\exists\; a_{1} , a_{2}\;,\;.\;.\;.\;,\;a_{n}\;\in$ F such that\\
 $w \in a_{1}\ast\alpha_{1}\;\#\;a_{2}\ast\alpha_{2}\;\#\;.\;.\;.\;\#\;a_{n}\ast\alpha_{n}$\\
Therefore $a\ast w\;\subseteq\;a\ast(a_{1}\ast\alpha_{1}\;\#\;a_{2}\ast\alpha_{2}\;\#\;.\;.\;.\;\#\;a_{n}\ast\alpha_{n})$\\
\smallskip\hspace{2.8cm}    $\subseteq\;a\ast(a_{1}\ast\alpha_{1})\;\#\;a\ast(a_{2}\ast\alpha_{2})\;\#\;.\;.\;.\;\#\;a\ast(a_{n}\ast\alpha_{n})$\\
\smallskip\hspace{2.8cm}    =$\;(a. a_{1})\ast\alpha_{1}\;\#\;(a. a_{2})\ast\alpha_{2}\;\#\;.\;.\;.\;\#(a. a_{n})\ast\alpha_{n}\;\subseteq$ W\\
Therefore $a\ast w\;\subseteq\;W\;,\;\forall\;a\in\;F$ and $\forall\;w\;\in$ W.\\
Hence W is a hypersubspace of V.\\
Next, since $0 , 1\, \in$ F. Let i $\in\;\{$1 , 2 , 3 , .  .  . , n$\}$ \\
Therefore \\$0\ast\alpha_{1}\;\#\; 0\ast\alpha_{2}\;\#\;.\;.\;.\;\#\;0\ast\alpha_{i-1}\;\#\;1\ast\alpha_{i}\;\#\;0\ast\alpha_{i+1}\;\#\;.\;.\;.\;\#\;0\ast\alpha_{n}\;\subseteq$ W\\
$\Rightarrow\;\theta\;\#\;\theta\;\#\;.\;.\;.\;\#\;\theta\;\#\;\alpha_{i}\;\#\;\theta\;.\;.\;.\;\#\;\theta\;\subseteq$ W\\
$\Rightarrow\;\alpha_{i}\;\in$ W , $[$ As $\theta\;\#\;\theta\;=\;\{\theta\}$ and $\theta\;\#\;\alpha_{i}\;=\;\alpha_{i}\;\#\;\theta\;=\;\alpha_{i}],$\\
therefore $\alpha_{i}\;\in$ W , for all i = 1 , 2 , .   .   . , n.\\
Let P be a hypersubspace of V containing $\alpha_{1}\;,\;\alpha_{2}\;,\;.\;.\;.\;,\;\alpha_{n}$.\\
Let $\alpha\in$ W , then there exists $a_{1},\;a_{2},\;.\;.\;.,\;a_{n}\;\in$ F such that \\
$\alpha\;\in\;a_{1}\ast\alpha_{1}\;\#\;a_{2}\ast\alpha_{2}\;\#\;.\;.\;.\;\#\;a_{n}\ast\alpha_{n}\;\subseteq$ P.\\
Therefore $\alpha\in$ P $\Rightarrow\;W\;\subseteq$ P\\
Hence W is the smallest hyperspace of V containing $\alpha_{1},\;\alpha_{2},\;.\;.\;.\;,\;\alpha_{n}$.

\begin{note} The linear combination of a null vector is the set $\{\theta\}$.
\end{note}

\section{\textbf{Linear Dependence and Linear Independence}}
\begin{definition} Let V be a hypervector space over a hyperfield F and $S =\{\alpha_{1} ,\alpha_{2} ,\;.\;.\;.,\;\alpha_{n}\}$ be a finite subset of V. Then S is said to be \textbf{linearly dependent} if there exist the scalars $a_{1} , a_{2} ,\;.\;.\;.\;,\;a_{n}\;\in$ F $($not all zero$)$ such that \\
$\theta \in a_{1}\ast\alpha_{1}\;\#\;a_{2}\ast\alpha_{2}\;\#\;.\;.\;.\;\#\;a_{n}\ast\alpha_{n}.$\\
Otherwise S is said to be \textbf{linearly independent}. i.e, if S is linearly independent and $\theta\in a_{1}\ast\alpha_{1}\;\#\;a_{2}\ast\alpha_{2}\;\#\;.\;.\;.\;\#\;a_{n}\ast\alpha_{n}.$\\
Then $a_{1}=a_{2}=.\;.\;.\;=a_{n}=0 $
\end{definition}
\begin{definition}Let V be a hypervector space over a hyperfield F and S $\subseteq$ V. Then S is said to be linearly dependent if S contains a finite subset which is linearly dependent. Otherwise S is linearly independent.
\end{definition}
\begin{result} Any singleton set of non-null vector of a hypervector space V is linearly independent.
\end{result}
\textbf{Proof}: Let $\alpha$ be a non-null vector of a hypervector space V. If possible let $\theta\;\in\;a\ast\alpha$ for some $a\in F$. We now show that $a=0$.\\
If $a \neq 0$, then $a^{-1}\in F$.\\
Now, since  $\theta \in a\ast \alpha\Rightarrow a^{-1}\ast\theta\subseteq a^{-1}\ast(a\ast\alpha)$\\
$\Rightarrow \theta\in a^{-1}\ast(a\ast\alpha)$    [ As $a^{-1}\ast\theta=\theta$]\\
$\Rightarrow \theta\in(a^{-1}.a)\ast\alpha$\\
$\Rightarrow \theta\in 1_{F}\ast\alpha$\\
$\Rightarrow\theta\in\{\alpha\}$
$\Rightarrow\alpha=\theta$, which is a contradiction.\\
Hence $a=0$. This completes the proof.

\begin{result} Any set of vectors containing the null vector is always linearly dependent.
\end{result}
\textbf{Proof}: Obvious.

\begin{definition}Let V be a hypervector space over a hyperfield F. Then the vector $\alpha\in$ V is said to be a \textbf{linear combination} of the vectors        $\alpha_{1},\alpha_{2}, .\; .\; . , \alpha_{n}\in$ V if there exist $a_{1},a_{2}, . \;.\; . ,a_{n}\in $F such that \\
$\alpha\in a_{1}\ast\alpha_{1}\;\#\;a_{2}\ast\alpha_{2}\;\#\; .\;.\; .\;\#\;a_{n}\ast\alpha_{n}$.
\end{definition}
\begin{theorem}\label{1} Let V be a hypervector space over a hyperfield F and S = $\{\alpha_{1},\alpha_{2}, .\; .\; . , \alpha_{n}\}$ be a subset of V. Then S is linearly dependent if and only if at least one of S can be expressed as a linear combination of the remaining other members of S.
\end{theorem}
\textbf{Proof}: Let S is linearly dependent. Then there exist $a_{1},a_{2},.\; .\; . ,a_{n}($ not all zero $)\in$ F such that \\\smallskip\hspace{1.5cm}
$\theta \in a_{1}\ast\alpha_{1}\;\#\;a_{2}\ast\alpha_{2}\;\#\;.\;.\;.\;\#\;a_{n}\ast\alpha_{n},$
\hspace{1.5cm}
$\cdots\hspace{1.5cm}(i)$\\
since $\#$ is commutative. Therefore without loss of generality we assume that $a_{1}\neq$ 0 . Then $a_{1}^{-1}\in$ F. Therefore from $(i)$ we get\\
$a_{1}^{-1}\ast\theta\subseteq a_{1}^{-1}\ast( a_{1}\ast\alpha_{1}\;\#\;a_{2}\ast\alpha_{2}\;\#\;.\;.\;.\;\#\;a_{n}\ast\alpha_{n})$\\
$\Rightarrow\;\{\theta\}\subseteq \{a_{1}^{-1}\ast( a_{1}\ast \alpha_{1})\}$ $\#\;\{a_{1}^{-1}\ast(a_{2}\ast\alpha_{2})\}\;\#\;.\;.\;.\;\#\;\{a_{1}^{-1}\ast(a_{n}\ast\alpha_{n})\}$\\
 \smallskip\hspace{1cm}=$(a_{1}^{-1}.a_{1})\ast\alpha_{1}\;\#\;(a_{1}^{-1}.a_{2})\ast\alpha_{2}\;\#\;.\; .\; .\#\;(a_{1}^{-1}.a_{n})\ast\alpha_{n}$\\
 \smallskip\hspace{1cm}=$1_{F}\ast\alpha_{1}\;\#\;(a_{1}^{-1}.a_{2})\ast\alpha_{2}\;\#\;.\; . \;.\#\;(a_{1}^{-1}.a_{n})\ast\alpha_{n}$\\
 \smallskip\hspace{1cm}=$\alpha_{1}\;\#\;(a_{1}^{-1}.a_{2})\ast\alpha_{2}\;\#\;. \;. \;.\#\;(a_{1}^{-1}.a_{n})\ast\alpha_{n}$\\
\smallskip\hspace{1cm}=$\alpha_{1}\;\#\;((a_{1}^{-1}.a_{2})\ast\alpha_{2}\;\#\;.\; .\; .\#\;(a_{1}^{-1}.a_{n})\ast\alpha_{n}).$\\
Then $\exists$ an element $\beta\in (a_{1}^{-1}.a_{2})\ast\alpha_{2}\;\#\;.\; . \;.\#\;(a_{1}^{-1}.a_{n})\ast\alpha_{n}$ such that $\theta\in\alpha_{1}\;\#\;\beta=\beta\;\#\;\alpha_{1}.$\\
This implies that $\alpha_{1}=-\beta$. So  $\alpha_{1}\,\in\,(-1_{F})\ast\beta.$
Therefore\\ $\alpha_{1}\,\in\,(-1_{F})\ast\beta\,\subseteq\,(-1_{F})\ast((a_{1}^{-1}.a_{2})\ast\alpha_{2}\;\#\;.\; .\; .\#\;(a_{1}^{-1}.a_{n})\ast\alpha_{n})$\\
\smallskip\hspace{3cm}$\subseteq(-1_{F})\ast((a_{1}^{-1}.a_{2})\ast\alpha_{2})\;\#\;.\; .\; .\#\;(-1_{F})\ast((a_{1}^{-1}.a_{n})\ast\alpha_{n})$\\
\smallskip\hspace{3cm}$=((-1_{F}).(a_{1}^{-1}.a_{2}))\ast\alpha_{2})\;\#\;.\; . \;.\#\;((-1_{F}).(a_{1}^{-1}.a_{n}))\ast\alpha_{n})$\\
\smallskip\hspace{3cm}$=(-a_{1}^{-1}.a_{2})\ast\alpha_{2}\;\#\;.\; . \;.\#\;(-a_{1}^{-1}.a_{n})\ast\alpha_{n}$\\
i.e $\;\alpha_{1}\in(-a_{1}^{-1}.a_{2})\ast\alpha_{2}\;\#\;.\; . \;.\#\;(-a_{1}^{-1}.a_{n})\ast\alpha_{n}$\\
$\Rightarrow\;\alpha_{1}\in\;HL(\alpha_{2},\alpha_{2}, .\; .\; . ,\alpha_{n})$\\
This completes the proof of the necessary part of the theorem.\\
Converse part\\
Without loss of generality we assume that $\alpha_{1}\in Hl(\alpha_{2},\alpha_{3}, .\; . \;. ,\alpha_{n}).$\\
Then $\exists\; a_{2},a_{3}, .\; .\; . ,a_{n}\in$ F such that \\
$\alpha_{1}\in a_{2}\ast\alpha_{2}\;\#\;a_{3}\ast\alpha_{3}\;\#\; .\; .\; .\;\# \;a_{n}\ast\alpha_{n}$\\
$\Rightarrow-\alpha_{1}\;\#\;\alpha_{1}\,\subseteq\, -\alpha_{1}\;\#\;a_{2}\ast\alpha_{2}\;\#\;a_{3}\ast\alpha_{3}\;\#\; .\; . \;.\;\#\; a_{n}\ast\alpha_{n}$\\
$\Rightarrow\theta\in-\alpha_{1}\;\#\;a_{2}\ast\alpha_{2}\;\#\;a_{3}\ast\alpha_{3}\;\#\; . \;. \;.\;\#\; a_{n}\ast\alpha_{n}$ $[$ As $\theta\in-\alpha_{1}\;\#\;\alpha_{1}]$\\
$\Rightarrow\;\theta\in(-1_{F})\ast\alpha_{1}\;\#\;a_{2}\ast\alpha_{2}\;\#\;a_{3}\ast\alpha_{3}\;\#\; .\; .\; .\;\#\; a_{n}\ast\alpha_{n}$\\
$\Rightarrow \theta\in a_{1}\ast\alpha_{1}\;\#\;a_{2}\ast\alpha_{2}\;\#\; . \;. \;.\;\#\; a_{n}\ast\alpha_{n}$ , where $a_{1}=-1_{F}(\neq 0).$\\
Therefore the set S=$\{\alpha_{1},\alpha_{2}, .\;.\; . ,\alpha_{n}\}$ is linearly dependent.

\begin{theorem}The non-zero vectors $\alpha_{1},\alpha_{2}, . . . ,\alpha_{n}$ of a hypervector space are linearly dependent if and only if one of them, say $\alpha_{i}$ , is a linear combination of the previous vectors , i.e $\alpha_{i}\in a_{1}\ast\alpha_{1}\;\#\;a_{2}\ast\alpha_{2}\;\#\; .\; .\; .\;\# \;a_{i-1}\ast\alpha_{i-1}$ , for some $a_{1}, a_{2}, . . . ,a_{i-1}\in$ F.
\end{theorem}
\textbf{Proof}:First we suppose that $\alpha_{1},\alpha_{2}, . . . ,\alpha_{n}$ are linearly dependent. Then we get a set of scalars $a_{1},a_{2}, . . . , a_{n}($not all zero$)\in$ F such that \\
\smallskip\hspace{1.5cm}$\theta\in a_{1}\ast\alpha_{1}\;\#\;a_{2}\ast\alpha_{2}\;\#\; . \;. \;.\;\#\; a_{n}\ast\alpha_{n}$ \hspace{2cm} $\cdots\hspace{1.5cm}(i)$\\
Let k be the largest integer such that $a_{k}\neq$ 0, here we see that k$\neq $1.\\
 If k=1 , then we see that $\alpha_{1}=\theta$, which contradicts the fact that $\alpha_{1},\alpha_{2}, . . . , \alpha_{n}$ are non-zero vectors.\\
Since k be the largest integer such that $a_{k}\neq$ 0, it follows that $a_{i}$=0 for all i, k$<i\leq$ n.\\
hence $(i)$ reduces to the following form \\
$\theta\in a_{1}\ast\alpha_{1}\;\#\;a_{2}\ast\alpha_{2}\;\#\; . \;. \;.\;\#\; a_{k}\ast\alpha_{k}$\\
$\Rightarrow\;\alpha_{k}\in(-a_{k}^{-1}. a_{1})\ast \alpha_{1}\;\#\;(-a_{k}^{-1}. a_{2})\ast\alpha_{2}\;\#\;.\;.\;.\;\#\;(-a_{k}^{-1} . a_{k-1})\ast\alpha_{k-1}$ , $[$ By the procedure of the proof of the Theorem \ref{1}$]$ \\
conversely , we suppose that $\alpha_{i}\in a_{1}\ast\alpha_{1}\;\#\;a_{2}\ast\alpha_{2}\;\#\; . \;. \;.\;\#\; a_{i-1}\ast\alpha_{i-1}$ , for some scalars $a_{1},a_{2}, . . . , a_{i-1}\in$ F.\\
Therefore the set $\{\alpha_{1},\alpha_{2}, . . . , \alpha_{i-1},\alpha_{i}\}$ of vectors is linearly dependent.\\
So $\{\alpha_{1},\alpha_{2}, . . . , \alpha_{i-1},\alpha_{i},\alpha_{i+1}, . . . ,\alpha_{n}\}$ of vectors is also linearly dependent.\\
This completes the proof.

\begin{theorem}$($Deletion theorem$):$ Let V be a strongly left distributive hypervector space over the hyperfield F and V be generated by a linearly dependent set $\{\alpha_{1},\alpha_{2},. . .,\alpha_{n}\}\subseteq$ V. Then V can also be generated by a suitable proper subset of $\{\alpha_{1},\alpha_{2},. . .,\alpha_{n}\}$.
\end{theorem}
\textbf{Proof}: Since V is generated by $\alpha_{1},\alpha_{2},. . .,\alpha_{n}$ , therefore we have V=HL$(\alpha_{1},\alpha_{2},. . .,\alpha_{n}).$\\
Again, since $\{\alpha_{1},\alpha_{2},. . .,\alpha_{n}\}$ is linearly dependent . It follows that one of the vectors $\alpha_{1},\alpha_{2},. . .,\alpha_{n}$ , say $\alpha_{i},$ can be expressed as a linear combination of the remaining others.\\
Then we get the scalars $c_{1},c_{2},. . .,c_{i-1},c_{i+1}, . . . ,c_{n}\in$ F such that \\
$\alpha_{i}\in c_{1}\ast\alpha_{1}\;\#\;c_{2}\ast\alpha_{2}\;\#\;.\;.\;.\;\#\;c_{i-1}\ast\alpha_{i-1}\;\#\;c_{i+1}\ast\alpha_{i+1}\;\#\;.\;.\;.\;\# c_{n}\ast\alpha_{n}$\\
We now show that HL$(\alpha_{1},\alpha_{2},...,\alpha_{i-1},\alpha_{i+1},...,\alpha_{n})$=HL$(\alpha_{1},\alpha_{2},...,\alpha_{i-1},\alpha_{i},\alpha_{i+1},...,\alpha_{n}).$\\
It is obvious that\\
HL$(\alpha_{1},\alpha_{2},...,\alpha_{i-1},\alpha_{i+1},...,\alpha_{n})\subseteq$ HL$(\alpha_{1},\alpha_{2},...,\alpha_{i},...,\alpha_{n})$=V.\\
Let $\alpha\in$ V, then there exist the scalars $a_{1},a_{2},...,a_{n}$ such that \\
$\alpha \in a_{1}\ast\alpha_{1}\;\#\;a_{2}\ast\alpha_{2}\;\#\; . \;. \;.\;\#\; a_{i}\ast\alpha_{i} \;\#\;.\;.\;.\;\#\;a_{n}\ast\alpha_{n}$\\\\
$\Rightarrow\;\alpha \in a_{1}\ast\alpha_{1}\;\#\;a_{2}\ast\alpha_{2}\;\#\; . \;. \;.\;\#\;a_{i-1}\ast\alpha_{i-1}\;\#\; a_{i}\ast (c_{1}\ast\alpha_{1}\;\#\;c_{2}\ast\alpha_{2}\;\#\;.\;.\;.\;\#\;c_{i-1}\ast\alpha_{i-1}\;\#\;c_{i+1}\ast\alpha_{i+1}\;\#\;.\;.\;.\;\# c_{n}\ast\alpha_{n} )\;\#\;a_{i+1}\ast\alpha_{i+1}\#\;.\;.\;.\;\#\;a_{n}\ast\alpha_{n}$\\\\
$\Rightarrow\;\alpha \in a_{1}\ast\alpha_{1}\;\#\;a_{2}\ast\alpha_{2}\;\#\; . \;. \;.\;\#\;a_{i-1}\ast\alpha_{i-1}\;\#\; (a_{i}. c_{1})\ast\alpha_{1}\;\#\;(a_{i}.c_{2})\ast\alpha_{2}\;\#\;.\;.\;.\;\#\;(a_{i}.c_{i-1})\ast\alpha_{i-1}\;\#\;(a_{i}.c_{i+1})\ast\alpha_{i+1}\;\#\;.\;.\;.\;\# (a_{i}.c_{n})\ast\alpha_{n} )\;\#\;a_{i+1}\ast\alpha_{i+1}\#\;.\;.\;.\;\#\;a_{n}\ast\alpha_{n}$\\\\
$\Rightarrow\;\alpha\in(a_{1}\oplus a_{i}. c_{1})\ast\alpha_{1}\;\#\;(a_{2}\oplus a_{i}.c_{2})\ast\alpha_{2}\;\#\;.\;.\;.\;\#\;(a_{i-1}\oplus a_{i}.c_{i-1})\ast\alpha_{i-1}\;\#\;(a_{i+1}\oplus a_{i}.c_{i+1})\ast\alpha_{i+1}\;\#\;.\;.\;.\;\# (a_{n}\oplus a_{i}.c_{n})\ast\alpha_{n}.$\\
Hence $\alpha\in$ HL$(\alpha_{1},\alpha_{2},...,\alpha_{i-1},\alpha_{i+1},...,\alpha_{n}).$\\
Therefore V=HL$(\alpha_{1},\alpha_{2},...,\alpha_{i-1},\alpha_{i+1},...,\alpha_{n}).$\\
This completes the proof.

\begin{theorem}If S=$\{\alpha_{1},\alpha_{2},...,\alpha_{n}\}$ is a linearly independent set of generators of a hypervector space V, then no proper subset of S can be a spanning set of V.
\end{theorem}
\textbf{Proof}: Obvious.

\begin{theorem}Suppose $\{\alpha_{1},\alpha_{2},...,\alpha_{n}\}$ generates a hypervector space V. If $\{\beta_{1},\beta_{2},...,\beta_{m}\}$ is linearly independent , then m$\leq$ n and V is generated by a set of the form $\{\beta_{1},\beta_{2},...,\beta_{m},\alpha_{i_{1}},\alpha_{i_{2}},...,\alpha_{i_{n-m}}\}$.
\end{theorem}
\textbf{Proof}: Obvious.

\section{\textbf{Basis or Hamel Basis}}
\begin{definition}Let V be a hypervector space over the hyperfield F and S be a subset of V. S is said to be a \textbf{basis, or Hamel basis} if\\
$(\,i\,)\;$ S is linearly independent.\\
$(\,ii\,)$ Every elements of V can be expressed as a finite linear combination of a few elements of S.\\\\
   If S is a basis of the hypervector space V and S is finite, then the hypervector space V is said to be a finite dimensional hypervector space and the  number of elements in S is called the dimension of the hypervector space V. Usually the dimension of V is denoted by dim$(V)$.\\
  Again if S is infinite Set , then V is said to be an infinite dimensional hypervector space.
\end{definition}
\begin{theorem} If S=$\{\alpha_{1},\alpha_{2},...,\alpha_{n}\}$ is a basis of a finite dimensional strongly left distributive hypervector space V over a hyperfield F, then every non-null vector $\alpha\in$ V has a unique representation.
\end{theorem}
\textbf{Proof}: Since S is a basis of V and $\alpha\in$ V, there exist $a_{1},a_{2},...,a_{n}\in$ F such that \\
\smallskip\hspace{1.5cm}$\alpha\in a_{1}\ast\alpha_{1}\;\#\;a_{2}\ast\alpha_{2}\;\#\;.\;.\;.\;\#\;a_{n}\ast\alpha_{n}.
$\hspace{2cm}$\cdots\hspace{1.5cm}(i)$\\
If possible let $\alpha\in b_{1}\ast\alpha_{1}\;\#\;b_{2}\ast\alpha_{2}\;\#\;.\;.\;.\;\#\;b_{n}\ast\alpha_{n}$  ,  for some $b_{1},b_{2},...,b_{n}\in$ F.\\
Therefore $-\alpha\,\in\,(-1_{F}\ast\alpha)\,\subseteq\,(-1_{F})\ast(b_{1}\ast\alpha_{1}\;\#\;b_{2}\ast\alpha_{2}\;\#\;.\;.\;.\;\#\;b_{n}\ast\alpha_{n})$\\
 $\Rightarrow\,-\alpha\in((-1_{F})\ast(b_{1}\ast\alpha_{1}))\;\#\;((-1_{F})\ast(b_{2}\ast\alpha_{2}))\#\;.\;.\;.\;\#\;((-1_{F})\ast(b_{n}\ast\alpha_{n}))$\\  \smallskip\hspace{.81cm}$=((-1_{F}).b_{1})\ast\alpha_{1}\;\#\;((-1_{F}).b_{2})\ast\alpha_{2}\#\;.\;.\;.\;\#\;((-1_{F}).b_{n})\ast\alpha_{n}$\\
\smallskip\hspace{.81cm}$=(-b_{1})\ast\alpha_{1}\;\#\;(-b_{2})\ast\alpha_{2}\#\;.\;.\;.\;\#\;(-b_{n})\ast\alpha_{n}$\\
Therefore $-\alpha\in(-b_{1})\ast\alpha_{1}\;\#\;(-b_{2})\ast\alpha_{2}\#\;.\;.\;.\;\#\;(-b_{n})\ast\alpha_{n}$ $\;\;\hspace{.4cm}\cdots\hspace{.6cm}(ii).$\\
From $(i)$ and $(ii)$ we get \\
$\alpha\,\#-\alpha\subseteq a_{1}\ast\alpha_{1}\;\#\;a_{2}\ast\alpha_{2}\;\#\;.\;.\;.\;\# a_{n}\ast\alpha_{n}\;\#\;(-b_{1})\ast\alpha_{1}\;\#\;(-b_{2})\ast\alpha_{2}\#\;.\;.\;.\;\#\;(-b_{n})\ast\alpha_{n}$\\
Therefore $\theta\in\alpha\,\#-\alpha\subseteq( a_{1}\oplus(-b_{1}))\ast\alpha_{1}\;\#\;(a_{2}\oplus(-b_{2}))\ast\alpha_{2}\;\#\;.\;.\;.\;\# (a_{n}\oplus(-b_{n}))\ast\alpha_{n}$.\\
Since $\{\alpha_{1},\alpha_{2},...,\alpha_{n}\}$ is a basis of V and $\theta\in(a_{1}\oplus(-b_{1}))\ast\alpha_{1}\;\#\;(a_{2}\oplus(-b_{2}))\ast\alpha_{2}\;\#\;.\;.\;.\;\# (a_{n}\oplus(-b_{n}))\ast\alpha_{n}.$\\
Then $0\in a_{i}\oplus (-b_{i})$ , for all i=1,2,...,n.\\
Again $(F\,,\,\oplus)$ is commutative.\\
Therefore $a_{i}=-(-b_{i})$ ,  for all i=1,2,...,n.\\
i.e $a_{i}=b_{i}$ ,  for all i=1,2,...,n.\\
This completes the proof.

\begin{theorem}$($Extension Theorem$): $ A linearly independent set of vectors in a finite dimension hypervector space V over a hyperfield F is either a basis, or it can be extended to a basis of V.
\end{theorem}
\textbf{Proof}:Let S=$\{\alpha_{1},\alpha_{2},...,\alpha_{n}\}$ be a linearly independent set of vectors in V.\\
Now HL$(S)$ being a smallest hypersubspace of V containing S. It therefore follows that HL$(S)\subseteq$ V.\\
If HL$(S)$=V, then S is a basis of V.\\
If HL$(S)$ is a proper hypersubspace of V, we show that S can be extended to a basis of V. \\
Let $\alpha\in$ V$\setminus$ HL$(S)$ and $S_{1}=\{\alpha_{1},\alpha_{2},...,\alpha_{n},\alpha\}.$\\
Now we consider the following expression :\\
\smallskip\hspace{1.5cm}$\theta\in c_{1}\ast\alpha_{1}\;\#\; c_{2}\ast\alpha_{2}\;\#\;.\;.\;.\;.\;\# c_{n}\ast\alpha_{n}\;\#\; c\ast\alpha$,\hspace{1 cm}$\cdots\hspace{1.5cm}(i)$
for some $c_{1},c_{2},...,c_{n},c\in F$
We now claim that c=0.\\
If c$\neq$ 0, then $c^{-1}\in$ F\\
Therefore $c^{-1}\ast\theta\subseteq c^{-1}\ast(c_{1}\ast\alpha_{1}\;\#\; c_{2}\ast\alpha_{2}\;\#\;.\;.\;.\;.\;\# c_{n}\ast\alpha_{n}\;\#\; c\ast\alpha)$\\
$\Rightarrow\;\theta\in (c^{-1}.c_{1})\ast\alpha_{1}\;\#\; (c^{-1}.c_{2})\ast\alpha_{2}\;\#\;.\;.\;.\;.\;\#(c^{-1}.c_{n})\ast\alpha_{n}\;\#\; (c^{-1}.c)\ast\alpha$\\
$\Rightarrow\;\theta\in (c^{-1}.c_{1})\ast\alpha_{1}\;\#\; (c^{-1}.c_{2})\ast\alpha_{2}\;\#\;.\;.\;.\;.\;\#(c^{-1}.c_{n})\ast\alpha_{n}\;\#\; 1_{F}\ast\alpha$\\
$\Rightarrow\;\theta\in (c^{-1}.c_{1})\ast\alpha_{1}\;\#\; (c^{-1}.c_{2})\ast\alpha_{2}\;\#\;.\;.\;.\;.\;\#(c^{-1}.c_{n})\ast\alpha_{n}\;\#\; \alpha.$\\
Then $\exists \,\beta\in(c^{-1}.c_{1})\ast\alpha_{1}\;\#\; (c^{-1}.c_{2})\ast\alpha_{2}\;\#\;.\;.\;.\;.\;\#(c^{-1}.c_{n})\ast\alpha_{n}$ such that $\theta\in\beta\;\#\;\alpha=\alpha\;\#\;\beta$\\
$\Rightarrow\alpha=-\beta\,\in\,(-1_{F}\ast\beta).$\\
Therefore $\alpha\in(-1_{F})\ast((c^{-1}.c_{1})\ast\alpha_{1}\;\#\; (c^{-1}.c_{2})\ast\alpha_{2}\;\#\;.\;.\;.\;.\;\#(c^{-1}.c_{n})\ast\alpha_{n})$\\ $\Rightarrow\alpha\in(-c^{-1}.c_{1})\ast\alpha_{1}\;\#\; (-c^{-1}.c_{2})\ast\alpha_{2}\;\#\;.\;.\;.\;.\;\#(-c^{-1}.c_{n})\ast\alpha_{n}$\\
$\Rightarrow \alpha\in HL(\alpha_{1},\alpha_{2},...,\alpha_{n})$\\
i.e $\alpha\in HL(S)$ , which is contradiction ,since $\alpha\in V\setminus HL(S).$
Hence we see that $c=0$ .\\
Substituting $c=0$ in $(i)$ we get \\
$\theta\in c_{1}\ast\alpha_{1}\;\#\; c_{2}\ast\alpha_{2}\;\#\;.\;.\;.\;.\;\# c_{n}\ast\alpha_{n}\;\#\; 0\ast\alpha$ \\
$\Rightarrow\theta\in c_{1}\ast\alpha_{1}\;\#\; c_{2}\ast\alpha_{2}\;\#\;.\;.\;.\;.\;\# c_{n}\ast\alpha_{n}\;\#\;\theta$\\
$\Rightarrow\theta\in c_{1}\ast\alpha_{1}\;\#\; c_{2}\ast\alpha_{2}\;\#\;.\;.\;.\;.\;\# c_{n}\ast\alpha_{n}$ \\
$\Rightarrow c_{1}=c_{2}=\;.\;.\;.\;=c_{n}=0,$ as $\alpha_{1},\alpha_{2},...,\alpha_{n}$ are linearly independent.\\
Thus we observe that $S_{1}$ is linearly independent.\\
If $HL(S_{1})$=V , then $S_{1}$ is a basis of V and the theorem is proved.\\
If $HL(S_{1})$ is a proper subspace of V , then again we can take a vector $\beta\in V\setminus HL(S_{1})$ and proceed as before .\\
Since V is a finite dimensional hypervector space , after a finite number of steps we come to a finite set of linearly independent vectors which generate the hypervector space V.\\
This completes the proof.

\begin{theorem} If $\{\alpha_{1},\alpha_{2},...,\alpha_{n}\}$ is a maximal linearly independent subset of a hypervector space V, then $\{\alpha_{1},\alpha_{2},...,\alpha_{n}\}$ is a basis of V.
\end{theorem}
\textbf{Proof}: Obvious.

\begin{theorem}Let U and W be two finite dimensional hypersubspace of a strongly left distributive hypervector space V over the hyperfield F, then U$\;\#\;$W is also a finite dimensional hypersubspace of V and\\
dim$(U\;\#\;W)=dim(U)+dim(W)-dim(U\cap W).$
\end{theorem}
\textbf{Proof}: Here we see that $U\cap W$ is a hypersubspace of both U and W. So $U\cap W$ is a finite dimensional hypersubspace of V. Let dim$(U\cap W)$=r , dim$(U)$=m and dim$(W)$=n. then we have r $\leq$ m , n. \\
Again let $\{\alpha_{1},\alpha_{2},...,\alpha_{r}\}$ be a basis of $U\cap W$.\\
Since $\{\alpha_{1},\alpha_{2},...,\alpha_{r}\}$ is a linearly independent set of vectors in U and r $\leq m=dim(U)$, then it follows that either  $\{\alpha_{1},\alpha_{2},...,\alpha_{r}\}$ is a basis of U or it can be extended to a basis for U.\\
Let  $\{\alpha_{1},\alpha_{2},...,\alpha_{r},\beta_{1},\beta_{2},...,\beta_{m-r}\}$ be a basis for U.\\
 By similar arguments we can suppose that  $\{\alpha_{1},\alpha_{2},...,\alpha_{r},\gamma_{1},\gamma_{2},...,\gamma_{n-r}\}$ is a basis of W.\\
Let S=$\{\alpha_{1},\alpha_{2},...,\alpha_{r},\beta_{1},\beta_{2},...,\beta_{m-r},\gamma_{1},\gamma_{2},...,\gamma_{n-r}\}.$\\
We now show that S is a basis of $U\;\#\;$W. \\
First we consider the following expression\\
$\theta\in a_{1}\ast\alpha_{1}\;\#\;a_{2}\ast\alpha_{2}\;\#\;.\;.\;.\;\#\;a_{r}\ast\alpha_{r}\;\#\;b_{1}\ast\beta_{1}\;\#\;b_{2}\ast\beta_{2}\;\#\;.\;.\;.\;\#\;
b_{m-r}\ast\beta_{m-r}\;\#\;c_{1}\ast\gamma_{1}\;\#\;c_{2}\ast\gamma_{2}\;\#\;.\;.\;.\;\#\;c_{n-r}\ast\gamma_{n-r}$,  \hspace{2cm} $\cdots\hspace{1.5cm}(i)$\\
 for some $a_{1},a_{2},...,a_{r},b_{1},b_{2},...,b_{m-r},c_{1},c_{2},...,c_{n-r}\in$ F.\\
Then $\exists $
 $\alpha\in a_{1}\ast\alpha_{1}\;\#\;a_{2}\ast\alpha_{2}\;\#\;.\;.\;.\;\#\;a_{r}\ast\alpha_{r}\;\#\;b_{1}\ast\beta_{1}\;\#\;\\
 \smallskip\hspace{4cm}b_{2}\ast\beta_{2}\;\#\;.\;.\;.\;\#\;
b_{m-r}\ast\beta_{m-r}$\smallskip\hspace{1cm} $\cdots\hspace{1.5cm}(ii)$\\
So $\theta\in\alpha\;\#\;c_{1}\ast\gamma_{1}\;\#\;c_{2}\ast\gamma_{2}\;\#\;.\;.\;.\;\#\;c_{n-r}\ast\gamma_{n-r}.$\\
Then $\exists\beta\in c_{1}\ast\gamma_{1}\;\#\;c_{2}\ast\gamma_{2}\;\#\;.\;.\;.\;\#\;c_{n-r}\ast\gamma_{n-r}$\\
So $\theta\in\alpha\;\#\;\beta=\beta\;\#\;\alpha,$    as $\#$ is commutative.\\
$\Rightarrow -\alpha=\beta$\\
$\Rightarrow -\alpha\in c_{1}\ast\gamma_{1}\;\#\;c_{2}\ast\gamma_{2}\;\#\;.\;.\;.\;\#\;c_{n-r}\ast\gamma_{n-r}$\hspace{2cm} $\cdots\hspace{1.5cm}(iii)$\\
$\Rightarrow (-1_{F})\ast(-\alpha)\subseteq(-1_{F})\ast (c_{1}\ast\gamma_{1}\;\#\;c_{2}\ast\gamma_{2}\;\#\;.\;.\;.\;\#\;c_{n-r}\ast\gamma_{n-r})$\\
$\Rightarrow -(-\alpha)\,\in\,((-1_{F}).c_{1})\ast\gamma_{1}\;\#\;((-1_{F}.c_{2})\ast\gamma_{2}\;\#\;.\;.\;.\;\#\;((-1_{F}).c_{n-r})\ast\gamma_{n-r}$\\
$\Rightarrow \alpha\in(-c_{1})\ast\gamma_{1}\;\#\;(-c_{2})\ast\gamma_{2}\;\#\;.\;.\;.\;\#\;(-c_{n-r})\ast\gamma_{n-r}$ \hspace{.5cm} $\cdots\hspace{1.5cm}(iv).$\\
From $(ii)$ and $(iv)$ we see that $\alpha\in U\cap W$.\\
Since $\alpha\in U\cap W$ and $\{\alpha_{1},\alpha_{2},...,\alpha_{r}\}$ is a basis of $U \cap W,$\\
we have $\alpha\in d_{1}\ast\alpha_{1}\;\#\;d_{2}\ast\alpha_{2}\;\;\#\;.\;.\;.\;\#\;d_{r}\ast\alpha_{r}$,  \hspace{2cm} $\cdots\hspace{1.5cm}(v)$\\
\smallskip\hspace{5cm} for some $d_{1},d_{2},...,d_{r}\in$ F. \\
Therefore from $(iii)$ and $(v)$ we get \\
$\alpha\;\#\;(-\alpha)\subseteq d_{1}\ast\alpha_{1}\;\#\;d_{2}\ast\alpha_{2}\;\#\;.\;.\;.\;\#\;d_{r}
\ast\alpha_{r}\;\#\;c_{1}\ast\gamma_{1}\;\#\;c_{2}\ast\gamma_{2}\;\#\;.\;.\;.\;\#\;c_{n-r}\ast\gamma_{n-r}$\\
$\Rightarrow\theta\in d_{1}\ast\alpha_{1}\;\#\;d_{2}\ast\alpha_{2}\;\#\;.\;.\;.\;\#\;d_{r}
\ast\alpha_{r}\;\#\;c_{1}\ast\gamma_{1}\;\#\;c_{2}\ast\gamma_{2}\;\#\;.\;.\;.\;\#\;c_{n-r}\ast\gamma_{n-r}$\\
$\Rightarrow d_{1}=d_{2}=...=d_{r}=c_{1}=c_{2}=...=c_{n-r}=0$, as $\{\alpha_{1},\alpha_{2},...,\alpha_{r},\gamma_{1},\gamma_{2},...,\gamma_{n-r}\}$ is a linearly independent set of vectors.\\
Then the expression $(i)$ reduces to the following\\
$\theta\in a_{1}\ast\alpha_{1}\;\#\;a_{2}\ast\alpha_{2}\;\#\;.\;.\;.\;\#\;a_{r}\ast\alpha_{r}\;\#\;b_{1}\ast
\beta_{1}\;\#\;b_{2}\ast\beta_{2}\;\#\;.\;.\;.\;\#\;b_{m-r}\ast\beta_{m-r}$\\
$\Rightarrow a_{1}=a_{2}=...=a_{r}=b_{1}=b_{2}=...=b_{m-r}=0$ ,  as$\{\alpha_{1},\alpha_{2},...,\alpha_{r},\beta_{1},\beta_{2},...,\beta_{m-r}\}$ is a linearly independent set of vectors.\\
Thus S is a linearly independent set of vectors V.\\
We now show that S generates $U\;\#\;W.$\\
Let $\alpha\in U\;\#\;W$\\
Then we have $\alpha\in\beta\;\#\;\gamma$ , for some $\beta\in$ U and for some $\gamma\in$ W.\\
Since $\beta\in$ U and $\{\alpha_{1},\alpha_{2},...,\alpha_{r},\beta_{1},\beta_{2},...,\beta_{m-r}\}$ is a basis of U.\\
Therefore $\beta\in a_{1}\ast\alpha_{1}\;\#\;a_{2}\ast\alpha_{2}\;\#\;.\;.\;.\;\#\;a_{r}\ast\alpha_{r}\;\#\;b_{1}\ast
\beta_{1}\;\#\;b_{2}\ast\beta_{2}\;\#\;.\;.\;.\;\#\;b_{m-r}\ast\beta_{m-r}$ , for some $a_{1},a_{2},...,a_{r},b_{1},b_{2},...,b_{m-r}\in$ F.\\
Again since $\gamma\in$ W and $\{\alpha_{1},\alpha_{2},...,\alpha_{r},\gamma_{1},\gamma_{2},...,\gamma_{n-r}\}$ is a basis of W.\\
Therefore $\gamma\in  c_{1}\ast\alpha_{1}\;\#\;c_{2}\ast\alpha_{2}\;\#\;.\;.\;.\;\#\;c_{r}\ast\alpha_{r}\;\#\;d_{1}\ast\gamma_{1}\;
\#\;d_{2}\ast\gamma_{2}\;\#\;.\;.\;.\;\#\;d_{n-r}\ast\gamma_{n-r}$ , for some $c_{1},c_{2},...,c_{r},d_{1},d_{2},...,d_{n-r}\in$ F.\\
Now $\alpha\in\beta\;\#\;\gamma.$\\
Therefore $\alpha\in a_{1}\ast\alpha_{1}\;\#\;a_{2}\ast\alpha_{2}\;\#\;.\;.\;.\;\#\;a_{r}\ast\alpha_{r}\;\#\;b_{1}\ast
\beta_{1}\;\#\;b_{2}\ast\beta_{2}\;\#\;.\;.\;.\;\#\;b_{m-r}\ast\beta_{m-r}\;\#\;c_{1}\ast\alpha_{1}\;\#\;c_{2}\ast
\alpha_{2}\;\#\;.\;.\;.\;\#\;c_{r}\ast\alpha_{r}\;\#\;d_{1}\ast\gamma_{1}\;\#\;d_{2}\ast\gamma_{2}\;\#\;.\;.\;.\;\#\;
d_{n-r}\ast\gamma_{n-r}$\\
$\Rightarrow\alpha\in (a_{1}\oplus c_{1})\ast\alpha_{1}\;\#\;(a_{2}\oplus c_{2})\ast\alpha_{2}\;\#\;.\;.\;.\;\#\;(a_{r}\oplus c_{r}\ast\alpha_{r}\;\#\;b_{1}\ast
\beta_{1}\;\#\;b_{2}\ast\beta_{2}\;\#\;.\;.\;.\;\#\;b_{m-r}\ast\beta_{m-r}\;\#\;d_{1}\ast\gamma_{1}\;
\#\;d_{2}\ast\gamma_{2}\;\#\;.\;.\;.\;\#\;d_{n-r}\ast\gamma_{n-r}.$\\
This shows that every vector of $U\;\#\;$W can be expressed as a linear combination of the vectors of S . Hence S is a basis of $U\;\#\;$W , which
proves that $U\;\#\;$W is a finite dimensional hypervector space of V.\\
Now dim$(U\;\#\;W)=|S|$\\
\smallskip\hspace{3.4cm}=r+m-r+n-r\\
\smallskip\hspace{3.4cm}=m+n-r\\
\smallskip\hspace{3.4cm}=dim$(U)+dim(W)-dim(U\cap W).$\\
This completes the proof.\\\\

\end{document}